\documentclass[11pt]{amsart}
\usepackage{amsmath,amssymb,mathtools}
\usepackage[letterpaper,hmargin=1.1in,vmargin=0.95in]{geometry}
\usepackage{microtype}
\usepackage[hidelinks]{hyperref}
\hypersetup{pdftitle={A proof of Bourin's crossed-Heinz conjecture}}
\newtheorem{theorem}{Theorem}[section]
\newtheorem{lemma}[theorem]{Lemma}
\newcommand{\C}{\mathbb C}
\newcommand{\Tr}{\operatorname{Tr}}
\newcommand{\uin}[1]{\lvert\!\lvert\!\lvert #1 \rvert\!\rvert\!\rvert}
\newcommand{\Uin}[1]{\bigl\lvert\!\bigl\lvert\!\bigl\lvert #1 \bigr\rvert\!\bigr\rvert\!\bigr\rvert}
\newcommand{\one}{\mathbf{1}}
\numberwithin{equation}{section}
\allowdisplaybreaks[2]
\title[Bourin's crossed-Heinz conjecture]{A proof of Bourin's crossed-Heinz conjecture}
\author[Yulong Tian]{Yulong Tian$^{*}$}
\address{Department of Mathematics, Northeast Forestry University, Harbin, P. R. China}
\email{tylmath@nefu.edu.cn}
\author{Jinli Xu}
\address{Department of Mathematics, Northeast Forestry University, Harbin, P. R. China}
\email{jclixv@qq.com}
\thanks{$^{*}$Corresponding author.}
\hypersetup{pdfauthor={Yulong Tian and Jinli Xu}}
\date{}
\begin{document}
\begin{abstract}
For positive semidefinite matrices $A,B$ of arbitrary order and every unitarily invariant norm $\uin{\cdot}$, we prove that $\uin{A^tB^{1-t}+B^tA^{1-t}}\le\uin{A+B}$ for all $0<t<1$, resolving Bourin's crossed-Heinz conjecture. We obtain this result by proving
\[
 \Uin{B^{(1-r)/2}A^rB^{(1-r)/2}+A^{(1-r)/2}B^rA^{(1-r)/2}}\le\Uin{A+B}
\]
for $1/2\le r<1$. We also establish the corresponding weak log-majorization in dimension two.
\end{abstract}
\maketitle

\section{Introduction}
Let $M_n(\C)$ denote the space of $n\times n$ complex matrices, and let $A,B\in M_n(\C)$ be positive semidefinite. For $0<t<1$, let $H_t(A,B)=A^tB^{1-t}+A^{1-t}B^t$ denote the Heinz sum of $A$ and $B$. The classical Heinz inequality \cite{BhatiaDavis} states that
\[
 \Uin{H_t(A,B)}\le\Uin{A+B}
\]
for every unitarily invariant norm $\uin{\cdot}$ on $M_n(\C)$. Reversing the order of the factors in the second summand, define
\[
 C_t(A,B)=A^tB^{1-t}+B^tA^{1-t}.
\]
Bourin \cite[Question~4]{Bourin} asked whether the crossed-Heinz inequality
\begin{equation}\label{eq:bourin}
 \Uin{C_t(A,B)}\le \Uin{A+B}
\end{equation}
also holds for every unitarily invariant norm. A natural approach would be to compare
$C_t(A,B)$ directly with $H_t(A,B)$. This approach was explored by Hayajneh and Kittaneh \cite{HK2013} and by Bhatia \cite{BhatiaTrace}. However, Bottazzi et al.~\cite[Example~3.1]{Bottazzi} showed that this comparison fails for $3\times3$ matrices in the operator norm.

Several inequalities related to Bourin's conjecture have been established in \cite{DHHK2023,HHK2017a,HHK2023,HHKL2022}. Hayajneh, Hayajneh and Kittaneh \cite{HHK2017b} proved \eqref{eq:bourin} for the trace norm and all $0<t<1$. For general unitarily invariant norms, the cases $t=1/4,3/4$ were settled in \cite{HHK2021}. Kittaneh and Ricard \cite{KR} subsequently extended the validity of \eqref{eq:bourin} to the entire interval $[1/4,3/4]$ by complex interpolation. More recently, Sababheh and Bourin \cite{SB} proved \eqref{eq:bourin} for $2\times2$ matrices and all $0<t<1$. In our recent work \cite{TX}, we established it for the operator and Frobenius norms in arbitrary dimension and for all $0<t<1$, and extended the operator-norm inequality to arbitrary $C^*$-algebras.

In this paper, we resolve Bourin's conjecture in full generality.

\begin{theorem}\label{thm:main}
Let $n\ge1$ and let $A,B\in M_n(\C)$ be positive semidefinite. If $0<t<1$, then
\[
 \Uin{A^tB^{1-t}+B^tA^{1-t}}\le\Uin{A+B}
\]
for every unitarily invariant norm on $M_n(\C)$.
\end{theorem}

Since $C_{1-t}(A,B)=C_t(A,B)^*$ and the central interval $[1/4,3/4]$ is already covered by \cite{KR}, it suffices to treat $3/4\le t<1$.

Following Kittaneh and Ricard \cite[Theorem~2.3]{KR}, Section~\ref{sec:truncation} uses a block factorization to reduce the problem to the inequality $\uin{K_r(A,B)}\le\uin{A+B}$ for every unitarily invariant norm, where
\[
 K_r(A,B)=B^{(1-r)/2}A^rB^{(1-r)/2}+A^{(1-r)/2}B^rA^{(1-r)/2},
 \qquad r=2t-1.
\]

Hayajneh, Hayajneh and Kittaneh \cite{HHK2023b} proved that
\[
 \Uin{K_r(A,B)}\le 2^{2r-1}\Uin{A+B},\qquad 1/2\le r<1,
\]
for every unitarily invariant norm. Theorem~\ref{thm:end-sandwich} improves the factor $2^{2r-1}$ to $1$.

For $1/2\le r\le5/6$, we prove the inequality by spectral truncation. Rather than comparing eigenvalue partial sums directly, we compare the threshold functions
\[
 E_h(Z)=\sum_{i=1}^n(\lambda_i(Z)-h)_+,
 \qquad h>0,
\]
where $\lambda_1(Z)\ge\cdots\ge\lambda_n(Z)$ are the eigenvalues of a Hermitian matrix $Z\in M_n(\C)$ and $x_+=\max\{x,0\}$. By the threshold characterization of weak majorization \cite[Theorem~6.3]{HiaiPetz}, it suffices to show that $E_h(K_r(A,B))\le E_h(A+B)$ for every $h>0$. For $u>0$, define
\[
 A_u=\min(A,uI),\qquad B_u=\min(B,uI),
\]
where the minima are defined by spectral calculus, and introduce the truncation paths
\[
 S(u)=A_u+B_u,\qquad K(u)=K_r(A_u,B_u).
\]
For each $h>0$, we show that the difference
\[
 F_h(u)=E_h(S(u))-E_h(K(u))
\]
is nonnegative at $3h$ and nondecreasing on $[3h,\infty)$. The initial comparison relies on a Hermitian congruence, while monotonicity follows by estimating the contributions of individual eigenvectors to the derivatives of the threshold functions. Since the truncated matrices eventually coincide with the original matrices, the desired inequality follows.

Section~\ref{sec:endpoints} treats the remaining range $5/6\le r<1$. We first derive a variational estimate in terms of the spectral deficit and the Jensen differences of fractional powers. Using resolvent integral representations, we then reduce the remaining matrix estimate to a scalar supporting-line inequality.

Section~\ref{sec:logmajorization} strengthens Theorem~\ref{thm:main} to weak log-majorization in dimension two, answering a question of Sababheh and Bourin \cite[Question~2.5]{SB} in this case.

\section{Reduction and spectral truncation}\label{sec:truncation}
For Hermitian matrices $X,Y$, we write $X\le Y$ if $Y-X$ is positive semidefinite. We write $\|\cdot\|$ for the Euclidean norm on vectors.

\subsection{Reduction to a positive matrix inequality}
Following the factorization in \cite[Theorem~2.3]{KR}, let
\[
 \mathcal R=\begin{pmatrix}A^{1/2}&B^{1/2}\end{pmatrix},\qquad
 \mathcal T=\begin{pmatrix}A^{t-1/2}B^{1-t}\\ B^{t-1/2}A^{1-t}\end{pmatrix}.
\]
Then $\mathcal R\mathcal R^*=A+B$ and $\mathcal R\mathcal T=C_t(A,B)$. By the Cauchy--Schwarz inequality for unitarily invariant norms \cite[Eq.~(IX.32)]{BhatiaMA}, applied to the rectangular factors by zero padding,
\begin{equation}\label{eq:end-factor}
 \begin{aligned}
 \Uin{C_t(A,B)}&\le\Uin{\mathcal R\mathcal R^*}^{1/2}\Uin{\mathcal T^*\mathcal T}^{1/2}\\
 &=\Uin{A+B}^{1/2}\Uin{B^{1-t}A^{2t-1}B^{1-t}+A^{1-t}B^{2t-1}A^{1-t}}^{1/2}.
 \end{aligned}
\end{equation}
Thus, setting $r=2t-1$, it suffices to establish the following theorem.

\begin{theorem}\label{thm:end-sandwich}
Let $A,B\in M_n(\C)$ be positive semidefinite. If $1/2\le r<1$, then
\[
 \Uin{B^{(1-r)/2}A^rB^{(1-r)/2}+A^{(1-r)/2}B^rA^{(1-r)/2}}\le\Uin{A+B}
\]
for every unitarily invariant norm on $M_n(\C)$.
\end{theorem}

\subsection{Threshold functions and the initial comparison}
Throughout the rest of this section, assume $1/2\le r\le5/6$, put $a=(1-r)/2$, and fix $h>0$.

\begin{lemma}\label{lem:end-initial}
We have $F_h(3h)\ge0$.
\end{lemma}
\begin{proof}
For a Hermitian matrix $T$, let $n_+(T)$ denote the number of its strictly positive eigenvalues, counted with multiplicity. By spectral calculus,
\[
 E_h(Z)=\sum_{i=1}^n\int_h^\infty\one_{\{v<\lambda_i(Z)\}}\,dv
 =\int_h^\infty n_+(Z-vI)\,dv,\qquad Z\ge0.
\]
Thus it suffices to compare $n_+(K(3h)-vI)$ and $n_+(S(3h)-vI)$ for every $v\ge h$. To reduce this comparison to a fixed threshold, take $X=2A_{3h}/v$ and $Y=2B_{3h}/v$. Since $v\ge h$, we have $0\le X,Y\le6I$. By the homogeneity of $K_r$, the required comparison will follow from a congruence inequality at threshold $2$.

For arbitrary $0\le X,Y\le6I$, denote $N=X^a+Y^a$. We claim that
\begin{equation}\label{eq:end-congruence}
 K_r(X,Y)-2I\le (N-I)(X+Y-2I)(N-I).
\end{equation}
The scalar inequality needed for this purpose is
\begin{equation}\label{eq:end-scalar-initial}
 x^{2-2a}(1-x^a)^2\le4x^a(1-x^a)(1-x^r),\qquad 0\le x\le6.
\end{equation}
The cases $x=0,1$ are immediate. For $0<x\le6$, $x\ne1$, the factors $1-x^{-a}$ and $1-x^{-r}$ have the same sign. Since $r\ge2a$, we have $|1-x^{-r}|\ge|1-x^{-2a}|=|1-x^{-a}|(1+x^{-a})$. Since $a\le1/4$ and $x^{-a}\ge6^{-a}\ge6^{-1/4}>1/2$, we obtain
\[
 x\frac{1-x^{-a}}{1-x^{-r}}\le\frac{x}{1+x^{-a}}
 \le\frac{6}{1+6^{-1/4}}<\frac{6}{1+1/2}=4.
\]
Using
\[
 \frac{1-x^{-a}}{1-x^{-r}}=x^{r-a}\frac{1-x^a}{1-x^r},
\]
multiplying by $x^a(1-x^a)(1-x^r)>0$ gives \eqref{eq:end-scalar-initial}, since $1+r=2-2a$.

Put
\[
 \begin{aligned}
 D_0&=\sum_{Z=X,Y}Z^{1-a}(I-Z^a),\\
 J_0&=\sum_{Z=X,Y}(I-Z^a)(I-Z^r).
 \end{aligned}
\]
For each $x\in[0,6]$, the scalar matrix
\[
 \begin{pmatrix}
 (1-x^a)(1-x^r)&x^{1-a}(1-x^a)/2\\
 x^{1-a}(1-x^a)/2&x^a
 \end{pmatrix}
\]
has nonnegative diagonal entries, and its determinant is nonnegative by \eqref{eq:end-scalar-initial}. Hence, by spectral calculus,
\[
 \begin{pmatrix}J_0&D_0/2\\D_0/2&N\end{pmatrix}
 =\sum_{Z=X,Y}
 \begin{pmatrix}
 (I-Z^a)(I-Z^r)&Z^{1-a}(I-Z^a)/2\\
 Z^{1-a}(I-Z^a)/2&Z^a
 \end{pmatrix}\ge0.
\]
We now express the difference in \eqref{eq:end-congruence} as a quadratic form in this block matrix. Since $r+2a=1$, we have
\[
 \begin{aligned}
 K_r(X,Y)&=\sum_{Z=X,Y}(N-Z^a)Z^r(N-Z^a)\\
 &=N(X^r+Y^r)N-N(X^{1-a}+Y^{1-a})\\
 &\quad-(X^{1-a}+Y^{1-a})N+X+Y.
 \end{aligned}
\]
Expanding $(N-I)(X+Y-2I)(N-I)$ and subtracting this identity yields
\[
 \begin{aligned}
 &2I+(N-I)(X+Y-2I)(N-I)-K_r(X,Y)\\
 &\quad=N(X+Y-X^r-Y^r-2I)N\\
 &\qquad+N(X^{1-a}+Y^{1-a}-X-Y)\\
 &\qquad+(X^{1-a}+Y^{1-a}-X-Y)N+4N\\
 &\quad=N(X+Y-X^r-Y^r-2I)N+ND_0+D_0N+4N.
 \end{aligned}
\]
On the other hand, expanding the block quadratic form gives
\[
 \begin{aligned}
 &\begin{pmatrix}N&2I-N\end{pmatrix}
 \begin{pmatrix}J_0&D_0/2\\D_0/2&N\end{pmatrix}
 \begin{pmatrix}N\\2I-N\end{pmatrix}\\
 &\quad=NJ_0N+ND_0+D_0N-ND_0N+(2I-N)N(2I-N)\\
 &\quad=N(J_0-D_0+N-4I)N+ND_0+D_0N+4N.
 \end{aligned}
\]
Since $J_0-D_0+N=X+Y-X^r-Y^r+2I$, the two expressions coincide. The block matrix is positive semidefinite, so this proves \eqref{eq:end-congruence}.

Returning to $X=2A_{3h}/v$ and $Y=2B_{3h}/v$, multiply \eqref{eq:end-congruence} by $v/2$ and use the homogeneity of $K_r$ to obtain
\[
 K(3h)-vI\le (N-I)(S(3h)-vI)(N-I),
\]
where $N=X^a+Y^a$ depends on $v$. The positive inertia index cannot increase under a decrease in the Loewner order or a Hermitian congruence, even when the congruence matrix $N-I$ is singular \cite[Corollaries~4.5.11 and 4.3.12]{HornJohnson}. Consequently,
\[
 n_+(K(3h)-vI)\le n_+(S(3h)-vI).
\]
Integrating over $v\in[h,\infty)$ proves the assertion.
\end{proof}

\subsection{Monotonicity along the truncation path}
\begin{lemma}\label{lem:end-propagation}
The function $F_h$ is nondecreasing on $[3h,\infty)$.
\end{lemma}
\begin{proof}
The paths $S(u)$ and $K(u)$ are locally Lipschitz for $u>0$, since the truncated matrices and their positive powers are given by scalar truncations in fixed eigenbases. The variational formula
\[
 E_h(Z)=\max_{0\le V\le I}\Tr(V(Z-hI)),\qquad Z\ge0,
\]
shows that $E_h$ is $n$-Lipschitz in the operator norm. Hence both threshold functions, and therefore $F_h$, are locally absolutely continuous. For $Z(u)=S(u)$ or $K(u)$, put $P_u=\one_{(h,\infty)}(Z(u))$. The variational formula gives
\[
 E_h(Z(v))-E_h(Z(u))\ge\Tr(P_u(Z(v)-Z(u))).
\]
At points where both $Z(u)$ and $E_h(Z(u))$ are differentiable, taking the right and left difference quotients yields
\[
 \frac{d}{du}E_h(Z(u))=\Tr(P_uZ'(u)).
\]
Fix $u\in(3h,\infty)\setminus(\operatorname{spec}(A)\cup\operatorname{spec}(B))$ at which both threshold functions are differentiable. Such points form a set of full measure in $[3h,\infty)$. Put
\[
 X=A_u/u,\quad Y=B_u/u,\quad
 Q=\one_{[0,h/u]}(X+Y),\quad P_K=\one_{(h,\infty)}(K(u)).
\]
Since $S(u)=u(X+Y)$, we have $\one_{(h,\infty)}(S(u))=I-Q$. We compare the contributions to the two derivatives one saturated spectral direction at a time. Choose orthonormal eigenbases of $A$ and $B$. A basis vector whose eigenvalue is below $u$ contributes zero to the derivative of every positive power of the corresponding truncated matrix. If $e$ is a unit basis vector of $A$ whose eigenvalue exceeds $u$, then $Xe=e$, and for every $\alpha>0$ its contribution to $(A_u^\alpha)'$ is $\alpha u^{\alpha-1}ee^*$. In particular, its contribution to $(E_h(S(u)))'$ is
\[
 \operatorname{Tr}((I-Q)ee^*)=1-\|Qe\|^2.
\]
We will show that the corresponding contribution to $(E_h(K(u)))'$ is no larger.

Since $Q$ commutes with $2I-X-Y$ and $Xe=e$,
\[
 (2I-X-Y)Qe=Q(2I-X-Y)e=Q(I-Y)e.
\]
On $\operatorname{ran}Q$, we have $2I-X-Y\ge(2-h/u)I\ge(5/3)I$ for $u\ge3h$. Hence, since $Q$ is a contraction,
\begin{equation}\label{eq:end-low}
 \|Qe\|^2\le\tfrac9{25}\|Q(I-Y)e\|^2
 \le\tfrac9{25}\langle e,(I-Y)^2e\rangle.
\end{equation}
In view of the preceding expression for the contribution to $(E_h(S(u)))'$, it is therefore enough to prove that the contribution of $e$ to $(E_h(K(u)))'$ is at most
\[
 1-\tfrac9{25}\langle e,(I-Y)^2e\rangle.
\]

Writing $A_u=uX$ and $B_u=uY$ and using $r+2a=1$ (so $u^{2a+r-1}=1$), the contribution of $e$ to $K'(u)$ is
\[
 rY^aee^*Y^a+a(ee^*Y^rX^a+X^aY^ree^*).
\]
Taking the trace against $P_K$, its contribution to $(E_h(K(u)))'=\operatorname{Tr}(P_KK'(u))$ is
\[
 r\|P_KY^ae\|^2+2a\operatorname{Re}\langle X^aP_Ke,Y^re\rangle.
\]
Since $X^ae=e$ and $X^a$ is a contraction, while $P_K$ is an orthogonal projection, we have
\[
 \|X^aP_Ke-e/2\|=\|X^a(P_K-I/2)e\|
 \le\|(P_K-I/2)e\|=1/2.
\]
Consequently, by the Cauchy--Schwarz inequality,
\[
 \begin{aligned}
 2\operatorname{Re}\langle X^aP_Ke,Y^re\rangle
 &=\langle e,Y^re\rangle
   +2\operatorname{Re}\langle X^aP_Ke-e/2,Y^re\rangle\\
 &\le\langle e,Y^re\rangle+\|Y^re\|.
 \end{aligned}
\]
Combining this with $\|P_KY^ae\|\le\|Y^ae\|$ and
$\|Y^re\|\le(1+\|Y^re\|^2)/2$, we obtain
\[
 \begin{aligned}
 &r\|P_KY^ae\|^2
   +2a\operatorname{Re}\langle X^aP_Ke,Y^re\rangle\\
 &\quad\le r\langle e,Y^{1-r}e\rangle
       +a\langle e,Y^re\rangle
       +\tfrac a2\langle e,(I+Y^{2r})e\rangle\\
 &\quad=1-\langle e,\psi_r(Y)e\rangle,
 \end{aligned}
\]
where
\[
 \psi_r(y)=r(1-y^{1-r})+a(1-y^r)+\tfrac a2(1-y^{2r}).
\]
Thus the desired bound for the contribution from $K'(u)$ will follow from the scalar estimate
\begin{equation}\label{eq:end-scalar-propagation}
 \psi_r(y)\ge\tfrac9{25}(1-y)^2,
 \qquad 0\le y\le1.
\end{equation}

To prove \eqref{eq:end-scalar-propagation}, the integral Cauchy--Schwarz inequality gives, for $0\le y\le1$,
\[
 \begin{aligned}
 (1-y)^2
 &=\left(\int_y^1x^{-r/2}x^{r/2}\,dx\right)^2\\
 &\le\int_y^1x^{-r}\,dx\int_y^1x^r\,dx
 =\frac{(1-y^{1-r})(1-y^{1+r})}{1-r^2}.
 \end{aligned}
\]
Since $1-y^{1+r}\le1$, this gives $(1-r^2)(1-y)^2\le1-y^{1-r}$. Using also
$1-y^r\ge r(1-y)$, $1-y^{2r}\ge1-y$, and $1-y\ge(1-y)^2$, we obtain
\[
 \begin{aligned}
 \psi_r(y)
 &\ge\left[r(1-r^2)+ar+\tfrac a2\right](1-y)^2\\
 &=(1-r)(r^2+3r/2+1/4)(1-y)^2\\
 &\ge\tfrac9{25}(1-y)^2.
 \end{aligned}
\]
Indeed, the coefficient $(1-r)(r^2+3r/2+1/4)$ is nonincreasing on $[1/2,5/6]$, since its derivative is $-3r^2-r+5/4\le0$. At $r=5/6$ it equals $79/216>9/25$, proving \eqref{eq:end-scalar-propagation}.

Since $0\le Y\le I$, spectral calculus gives
$\psi_r(Y)\ge\tfrac9{25}(I-Y)^2$. Together with \eqref{eq:end-low}, this gives
\[
 \begin{aligned}
 &r\|P_KY^ae\|^2
   +2a\operatorname{Re}\langle X^aP_Ke,Y^re\rangle\\
 &\quad\le1-\langle e,\psi_r(Y)e\rangle\\
 &\quad\le1-\tfrac9{25}\langle e,(I-Y)^2e\rangle\\
 &\quad\le1-\|Qe\|^2.
 \end{aligned}
\]
Thus, the summand corresponding to $e$ in $\Tr(P_KK'(u))$ does not exceed the corresponding summand in $\Tr((I-Q)S'(u))$.
Interchanging $A$ and $B$ and summing over the corresponding basis vectors gives
\[
 F_h'(u)\ge0\qquad\text{almost everywhere on }[3h,\infty).
\]
Local absolute continuity proves the claim.
\end{proof}

By Lemmas~\ref{lem:end-initial} and \ref{lem:end-propagation}, $F_h(u)\ge F_h(3h)\ge0$ for every $u\ge3h$. If $u\ge\max\{3h,\lambda_1(A),\lambda_1(B)\}$, then
\[
 E_h(K_r(A,B))\le E_h(A+B).
\]
Since $h>0$ is arbitrary, this proves Theorem~\ref{thm:end-sandwich} for $1/2\le r\le5/6$.

\section{A variational approach to the endpoint interval}\label{sec:endpoints}
Write $\|Z\|_{(k)}=\sum_{i=1}^k\sigma_i(Z)$ for the $k$th Ky Fan norm of $Z\in M_n(\C)$. In Subsections~\ref{sec:end-variational}--\ref{sec:end-scalar}, assume that $A,B$ are positive definite and $5/6\le r<1$. By the Ky Fan dominance theorem \cite[Theorem~IV.2.2]{BhatiaMA} and continuity, it remains to prove
\begin{equation}\label{eq:end-sandwich}
 \|K_r(A,B)\|_{(k)}\le \|A+B\|_{(k)},\qquad 1\le k\le n.
\end{equation}

\subsection{A variational estimate}\label{sec:end-variational}
Set
\[
 M=\frac{A^r+B^r}{2},\qquad
 \theta=\frac{1-r}{2r},\qquad p=1+2\theta=\frac1r.
\]
Since $5/6\le r<1$, we have $0<\theta\le1/10$ and $1<p\le6/5$. For $v\in\mathbb R$, write
\[
 M_v=\frac{A^{rv}+B^{rv}}2.
\]
Here the subscript in $M_v$ indicates the average of the $rv$th powers of $A$ and $B$, whereas $M^v$ denotes the $v$th power of their average $M$.
Define the Jensen differences associated with the power functions $x^p$, $x^{2\theta}$, and $x^\theta$ by
\[
 \begin{aligned}
 G&=M_p-M^p,\\
 J&=M^{2\theta}-M_{2\theta},\\
 D&=M^\theta-M_\theta,
 \end{aligned}
\]
with signs chosen so that each matrix is positive semidefinite. Indeed, since $1<p\le6/5<2$ and $0<\theta<2\theta\le1/5<1$, the function $x^p$ is operator convex \cite[Theorem~1.5.8]{BhatiaPDM}, while $x^\theta$ and $x^{2\theta}$ are operator concave on $(0,\infty)$ \cite[Theorems~V.2.5 and V.2.10]{BhatiaMA}; consequently, $G\ge0$, $J\ge0$, and $0\le D\le M^\theta$ (as $M_\theta>0$). Also define the scalar constant
\[
 b=\max_{0\le x\le1}(1-x)x^\theta
   =\frac{\theta^\theta}{(1+\theta)^{1+\theta}}.
\]
The choice of $b$ gives $(1-x)_+x^\theta\le b$ for all $x\ge0$, which will be used in Lemma~3.1 to control the lower spectral part.

Now fix $1\le k\le n$. By homogeneity, we may assume that the $k$th largest eigenvalue of $M$ is $1$. By Ky Fan's maximum principle, it suffices to show that $\Tr(PK_r(A,B))\le \|A+B\|_{(k)}$ for every rank-$k$ orthogonal projection $P$. Let $Q$ be a rank-$k$ orthogonal projection onto a subspace spanned by eigenvectors corresponding to the $k$ largest eigenvalues of $M$. Put $T=M^\theta$. Since $T^2=M^{2\theta}$ and $x\mapsto x^{2\theta}$ is increasing, $Q$ is also a maximizing projection for $T^2$. Thus
\[
 \Tr(QT^2)=\|T^2\|_{(k)}.
\]
Define the spectral deficit of $P$ by
\[
 \delta=\Tr((Q-P)T^2)=\|T^2\|_{(k)}-\Tr(PT^2)\ge0.
\]

\begin{lemma}\label{lem:end-variational}
Under the above normalization,
\begin{equation}\label{eq:end-variational}
 \frac{\|A+B\|_{(k)}-\Tr(PK_r(A,B))}2\ge\delta+2\Tr(P(J-bD)).
\end{equation}
\end{lemma}
\begin{proof}
Set $W=(M-I)_+$ and $L=(I-M)_+$. Since the $k$th largest eigenvalue of $M$ is $1$, we have $QW=W$ and $QL=0$, even if the eigenvalue $1$ has multiplicity greater than one. For a Hermitian matrix $Z$ define
\[
 \begin{aligned}
 \mathcal E(Z)&=\Tr((I-P)ZWZ)+\Tr(PZLZ)\\
 &=\|(I-P)ZW^{1/2}\|_{\mathrm F}^2+\|PZL^{1/2}\|_{\mathrm F}^2,
 \end{aligned}
\]
where $\|\cdot\|_{\mathrm F}$ denotes the Frobenius norm. By cyclicity of the trace and $W-L=M-I$,
\[
 \begin{aligned}
 \mathcal E(Z)&=\Tr(WZ^2)-\Tr(PZ(W-L)Z)\\
 &=\Tr((W+P)Z^2)-\Tr(PZMZ).
 \end{aligned}
\]
In particular, since $T$ commutes with $M$ and $W=Q(M-I)$,
\[
 \mathcal E(T)=\Tr((Q-P)(M-I)T^2).
\]
Expanding $M$ and using $rp=1$ and $r\theta=(1-r)/2$, we obtain
\[
 A^{(1-r)/2}MA^{(1-r)/2}
 +B^{(1-r)/2}MB^{(1-r)/2}
 =M_p+\tfrac12 K_r(A,B).
\]
Thus, applying the identity for $\mathcal E(Z)$ to $A^{(1-r)/2}$ and $B^{(1-r)/2}$ and using $\|A+B\|_{(k)}/2\ge\Tr(QM_p)$, we obtain
\[
 \begin{aligned}
 \frac{\|A+B\|_{(k)}-\Tr(PK_r(A,B))}2
 &\ge\Tr((Q+P)M_p)-2\Tr((W+P)M_{2\theta})\\
 &\quad+\mathcal E(A^{(1-r)/2})+\mathcal E(B^{(1-r)/2}).
 \end{aligned}
\]
To isolate the spectral deficit, use $M^p=MT^2$ and $W=Q(M-I)$ to write
\[
 \begin{aligned}
 &\Tr((Q+P)M^p)-2\Tr((W+P)T^2)\\
 &\quad=\Tr((Q+P)MT^2)-2\Tr((Q(M-I)+P)T^2)\\
 &\quad=\Tr((Q-P)(2I-M)T^2)=\delta-\mathcal E(T).
 \end{aligned}
\]
Since $\mathcal E$ is a quadratic form and
\[
 \frac{A^{(1-r)/2}+B^{(1-r)/2}}2=M_\theta=T-D,
\]
the pair $(A^{(1-r)/2},B^{(1-r)/2})$ and the pair $(T,T-2D)$ have the same midpoint. Applying the parallelogram identity to both pairs gives
\[
 \begin{aligned}
 \mathcal E(A^{(1-r)/2})+\mathcal E(B^{(1-r)/2})
 &=2\mathcal E(T-D)+\tfrac12\mathcal E(A^{(1-r)/2}-B^{(1-r)/2})\\
 &=\mathcal E(T)-2\mathcal E(D)+\mathcal E(T-2D)\\
 &\quad+\tfrac12\mathcal E(A^{(1-r)/2}-B^{(1-r)/2}).
 \end{aligned}
\]
Since $\mathcal E(Z)\ge0$ for every Hermitian $Z$, it follows that
\[
 \mathcal E(A^{(1-r)/2})+\mathcal E(B^{(1-r)/2})-\mathcal E(T)+2\mathcal E(D)\ge0.
\]
Substituting $M_p=M^p+G$ and $M_{2\theta}=T^2-J$, we obtain
\begin{equation}\label{eq:end-deficit}
 \begin{aligned}
 \frac{\|A+B\|_{(k)}-\Tr(PK_r(A,B))}2
 \ge{}&\delta+\Tr((Q+P)G)+2\Tr((W+P)J)-2\mathcal E(D).
 \end{aligned}
\end{equation}
We estimate the two summands of $\mathcal E(D)$ separately.

Write $D=T^{1/2}FT^{1/2}$, where $0\le F\le I$. In an eigenbasis of $M$, let $m_i$ be its eigenvalues and put $u_i=(m_i-1)_+m_i^\theta$, $v_i=m_i^\theta$. In this basis,
\[
 \Tr(WTD)=\sum_i u_iv_iF_{ii},\qquad
 \Tr(WD^2)=\sum_{i,j}u_iv_j|F_{ij}|^2.
\]
Using the symmetry $|F_{ij}|^2=|F_{ji}|^2$ to symmetrize the double sum, we obtain
\[
 \begin{aligned}
 \Tr(WTD)-\Tr(WD^2)
 &=\sum_i u_iv_i(F-F^2)_{ii}\\
 &\quad+\frac12\sum_{i,j}(u_i-u_j)(v_i-v_j)|F_{ij}|^2\ge0,
 \end{aligned}
\]
since $F^2\le F$ and $u_i,v_i$ are nondecreasing functions of $m_i$. Hence
\begin{equation}\label{eq:end-upper-trace}
 \Tr((I-P)DWD)\le\Tr(WD^2)\le\Tr(WTD).
\end{equation}
The scalar inequality
\[
 x^{1+2\theta}-x-2x^{2\theta}+2x^\theta
 =x^\theta(x^\theta-1)(x+x^{1-\theta}-2)\ge0
\]
for $x\ge0$ follows because the two factors in parentheses have the same sign. For a unit eigenvector $e$ of $M$ with eigenvalue $m>1$, we have
\[
 \langle e,(TD-J)e\rangle=\langle e,(M_{2\theta}-m^\theta M_\theta)e\rangle.
\]
Applying the scalar inequality to $A^r/m$ and $B^r/m$ by spectral calculus, taking the average and multiplying by $m^{2\theta}$ gives
\[
 m^{-1}M_p-m^{2\theta-1}M\ge2(M_{2\theta}-m^\theta M_\theta).
\]
Since $Me=me$ and $p=1+2\theta$, we have $m^{2\theta-1}Me=m^{-1}M^pe$. Taking the inner product against $e$ therefore yields
\[
 2\langle e,(TD-J)e\rangle\le m^{-1}\langle e,Ge\rangle.
\]
Multiplying by $m-1>0$ and using $(m-1)/m<1$ together with $G\ge0$ yields $2(m-1)\langle e,(TD-J)e\rangle\le\langle e,Ge\rangle$.
Summing over an orthonormal eigenbasis of $M$ restricted to $m>1$ gives $2\Tr(W(TD-J))\le\Tr(QG)$. Together with \eqref{eq:end-upper-trace}, this yields
\[
 2\Tr((I-P)DWD)\le\Tr(QG)+2\Tr(WJ).
\]
For the other summand, $(1-x)_+x^\theta\le b$ for all $x\ge0$ by the definition of $b$, so $L\le bT^{-1}$ and
\[
 DLD\le bDT^{-1}D=bT^{1/2}F^2T^{1/2}\le bD.
\]
Combining the two estimates gives
\[
 \begin{aligned}
 2\mathcal E(D)&=2\Tr((I-P)DWD)+2\Tr(PDLD)\\
 &\le\Tr(QG)+2\Tr(WJ)+2b\Tr(PD).
 \end{aligned}
\]
Substituting this bound into \eqref{eq:end-deficit} cancels the terms $\Tr(QG)$ and $2\Tr(WJ)$. Discarding the nonnegative term $\Tr(PG)$ then proves \eqref{eq:end-variational}.
\end{proof}

\subsection{Resolvent reduction and the matrix estimate}
By operator convexity of inversion \cite[Eq.~(1.33)]{BhatiaPDM}, for $u>0$,
\[
 R(u):=\frac{(A^r+uI)^{-1}+(B^r+uI)^{-1}}2-(M+uI)^{-1}\ge0.
\]
For $0<v<1$, the integral representation \cite[p.~22]{BhatiaPDM}
\[
 x^v=c_v\int_0^\infty u^{v-1}\frac{x}{x+u}\,du,\qquad c_v=\frac{\sin(\pi v)}{\pi},
\]
together with $x/(x+u)=1-u/(x+u)$ gives $M^v-M_v=c_v\int_0^\infty u^vR(u)\,du$. Define
\[
 w(u)=2(c_{2\theta}u^{2\theta}-bc_\theta u^\theta),\qquad
 g(x)=\int_0^\theta\frac{w(u)}{x+u}\,du\quad(x\ge0).
\]
For use below, put
\[
 M_g=\frac{g(A^r)+g(B^r)}2.
\]
Since $w(u)=-2bc_\theta u^\theta+O(u^{2\theta})$ as $u\downarrow0$, the function $g$ is well defined and continuous on $[0,\infty)$, and $w(u)<0$ for all sufficiently small $u>0$. Thus the positivity of $R(u)$ alone does not give a positive semidefinite lower bound for the full integral. On the other hand, since $c_{2\theta}/c_\theta=2\cos(\pi\theta)>1$ and $b<\theta^\theta$, we have $w(u)\ge0$ for $u\ge\theta$. As $R(u)\ge0$, the integrand is positive semidefinite on $[\theta,\infty)$, so discarding this tail in a lower bound gives
\begin{equation}\label{eq:end-resolvent}
 2(J-bD)=\int_0^\infty w(u)R(u)\,du
 \ge\int_0^\theta w(u)R(u)\,du=M_g-g(M).
\end{equation}
By \eqref{eq:end-variational} and \eqref{eq:end-resolvent}, it suffices to show that $\Tr(P(g(M)-M_g))\le\delta$. The following scalar estimate yields this bound.

\begin{lemma}\label{lem:end-envelope}
Let $0<\theta\le1/10$ and set $f(s)=g(s^{1/(2\theta)})$ for $s\ge0$. For every $\lambda\ge1$, we have $f'(\lambda)\le0$ and
\begin{equation}\label{eq:end-envelope}
 f(s)\ge f(\lambda)+f'(\lambda)(s-\lambda),\qquad s\ge0.
\end{equation}
Moreover, the function $s\mapsto s+f(s)$ is increasing on $[0,\infty)$.
\end{lemma}
We defer the proof of Lemma~\ref{lem:end-envelope} to Subsection~\ref{sec:end-scalar} and use it to derive the endpoint estimate.

Under the normalization of Subsection~\ref{sec:end-variational}, let $\lambda_1\ge\cdots\ge\lambda_k$ be the $k$ largest eigenvalues of $T^2$, and let $e_1,\ldots,e_k$ be an orthonormal eigenbasis of its compression to $\operatorname{ran}P$, ordered by decreasing eigenvalues. Then $\lambda_i\ge1$, and by Cauchy's interlacing theorem \cite[Corollary~III.1.5]{BhatiaMA},
\[
 \langle e_i,M_{2\theta}e_i\rangle\le\langle e_i,T^2e_i\rangle\le\lambda_i,\qquad 1\le i\le k.
\]

Applying \eqref{eq:end-envelope} to $A^{1-r}$ and $B^{1-r}$ by spectral calculus and taking the average, we obtain
\[
 \begin{aligned}
 \langle e_i,M_ge_i\rangle
 &\ge f(\lambda_i)+f'(\lambda_i)(\langle e_i,M_{2\theta}e_i\rangle-\lambda_i)\\
 &\ge f(\lambda_i),
 \end{aligned}
\]
since $f'(\lambda_i)\le0$. Since $g(M)=f(T^2)$, summing over $i$ gives
\[
 \Tr(PM_g)=\sum_{i=1}^k\langle e_i,M_ge_i\rangle\ge\sum_{i=1}^k f(\lambda_i)=\Tr(Qg(M)).
\]
On the other hand, since $s+f(s)$ is increasing, $Q$ is also a projection onto the top $k$ eigenvectors of $T^2+g(M)=T^2+f(T^2)$. By Ky Fan's maximum principle,
\[
 \Tr(P(T^2+g(M)))\le\Tr(Q(T^2+g(M))),
\]
so that $\Tr(Pg(M))-\Tr(Qg(M))\le\delta$. Together these two inequalities give
\[
 \Tr(P(g(M)-M_g))\le\delta.
\]
Combining \eqref{eq:end-variational} and \eqref{eq:end-resolvent}, we obtain
\[
 \frac{\|A+B\|_{(k)}-\Tr(PK_r(A,B))}2
 \ge\delta-\Tr(P(g(M)-M_g))\ge0.
\]
Taking the maximum over all rank-$k$ projections $P$ reduces \eqref{eq:end-sandwich} for $5/6\le r<1$ to Lemma~\ref{lem:end-envelope}, whose proof is given next.

\subsection{Proof of the scalar supporting-line estimate}\label{sec:end-scalar}
\begin{proof}[Proof of Lemma~\ref{lem:end-envelope}]
Substituting $u=\theta v$ gives
\[
 g(x)=\gamma F(x/\theta),\qquad
 F(z)=\int_0^1\frac{\alpha v^{2\theta}-v^\theta}{z+v}\,dv,
\]
where $\gamma=2bc_\theta\theta^\theta$ and $\alpha=c_{2\theta}\theta^\theta/(bc_\theta)=2\cos(\pi\theta)(1+\theta)^{1+\theta}$.

\smallskip\noindent\textit{Step 1: Global supporting lines.}\enspace By Bernoulli's inequality, $(1+\theta)^{1+\theta}\ge1+\theta+\theta^2$. Together with $\cos(\pi\theta)\ge1-5\theta^2$, this gives
\[
 \alpha\ge2(1-5\theta^2)(1+\theta+\theta^2)>2+\theta,
\]
since $8\theta+10\theta^2+10\theta^3\le91/100<1$ for $0<\theta\le1/10$.

To establish a global supporting line at $\lambda\ge1$, we will show that $f$ is concave on $[0,s_0]$ and convex on $[s_0,\infty)$ for some $s_0<1$, and that the tangent at $\lambda$ lies below $f(0)$ at $s=0$. The intercept of the tangent to $f$ at $s>0$ is $f(s)-sf'(s)$, whose derivative is $-sf''(s)$. It is therefore useful to study
\[
 \Phi(z)=2\theta F(z)-zF'(z),\qquad z>0.
\]
Indeed, by the chain rule,
\[
 f(s)-sf'(s)=\frac{\gamma}{2\theta}\Phi\left(\frac{s^{1/(2\theta)}}{\theta}\right),\qquad s>0.
\]
Since
\[
 \frac{d}{dv}\left(\frac{v^{1+\beta}}{z+v}\right)
 =v^\beta\left(\frac{\beta}{z+v}+\frac{z}{(z+v)^2}\right),
 \qquad \beta\in\{\theta,2\theta\},
\]
integrating over $v\in[0,1]$ gives
\[
 \begin{aligned}
 \Phi(z)&=\int_0^1(\alpha v^{2\theta}-v^\theta)
 \left(\frac{2\theta}{z+v}+\frac{z}{(z+v)^2}\right)\,dv\\
 &=\frac{\alpha-1}{z+1}-\theta\int_0^1\frac{v^\theta}{z+v}\,dv.
 \end{aligned}
\]
Differentiating, we obtain
\[
 -(z+1)^2\Phi'(z)=\alpha-1-\theta\int_0^1v^\theta\left(\frac{z+1}{z+v}\right)^2\,dv.
\]
The right-hand side is strictly increasing in $z$, since $(z+1)/(z+v)$ is strictly decreasing in $z$ for each $0<v<1$. As $z\downarrow0$, the integral increases to $\int_0^1 v^{\theta-2}\,dv=+\infty$ by monotone convergence. At $z=1$ the right-hand side is at least
\[
 \alpha-1-4\theta\int_0^1\frac{dv}{(1+v)^2}=\alpha-1-2\theta>0.
\]
Consequently, $-\Phi'(z)$ changes sign from negative to positive at a unique point $z_0\in(0,1)$. Differentiating the intercept identity above, we see that $f''(s)$ has the sign of $-\Phi'(s^{1/(2\theta)}/\theta)$. Thus, with $s_0=(\theta z_0)^{2\theta}<\theta^{2\theta}<1$, the function $f$ is concave on $[0,s_0]$ and convex on $[s_0,\infty)$, using continuity at $s=0$.

Fix $\lambda\ge1$ and put $z=\lambda^{1/(2\theta)}/\theta\ge1/\theta$. Since $\alpha>2+\theta$ and $F(0)=(\alpha-2)/(2\theta)$, we have
\[
 \begin{aligned}
 \frac{2\theta}{\gamma}\bigl(f(\lambda)-\lambda f'(\lambda)\bigr)
 &=\Phi(z)\le\frac{\alpha-1}{z+1}\\
 &\le\frac{\theta(\alpha-1)}{1+\theta}<\alpha-2=\frac{2\theta}{\gamma}f(0).
 \end{aligned}
\]
Since $s_0<1\le\lambda$, convexity gives \eqref{eq:end-envelope} for $s\ge s_0$. On $[0,s_0]$, the function $f(s)-f(\lambda)-f'(\lambda)(s-\lambda)$ is concave and nonnegative at both endpoints, so it is nonnegative throughout. This proves \eqref{eq:end-envelope}. Since $F(z)=O(1/z)$ as $z\to\infty$, we have $f(s)\to0$ as $s\to\infty$. Thus \eqref{eq:end-envelope} implies $f'(\lambda)\le0$.

\smallskip\noindent\textit{Step 2: Monotonicity of $s+f(s)$.}\enspace Applying Bernoulli's inequality to $(1+\theta)^{-\theta}$ gives
$(1+\theta)^{-1-\theta}=(1+\theta)^{-\theta}/(1+\theta)\ge(1-\theta^2)/(1+\theta)=1-\theta$.
Also, by concavity of sine, $\sin(\pi\theta/2)\ge\theta$, and hence
$\cos(\pi\theta)=1-2\sin^2(\pi\theta/2)\le1-2\theta^2$. Therefore
\[
 \frac{\alpha}{2}\le\frac{1-2\theta^2}{1-\theta}<\frac{1+2\theta}{1+\theta}.
\]

Writing $(\theta z)^{2\theta}=\theta^{2\theta}z^{2\theta}$, differentiating the integral representation of $z^{2\theta}$ and using $\theta^{2\theta}c_{2\theta}/\gamma=\alpha/2$, we obtain
\[
 \begin{aligned}
 \frac1\gamma\frac{d}{dz}\bigl[(\theta z)^{2\theta}+\gamma F(z)\bigr]
 &=\frac{\alpha}{2}\int_0^\infty\frac{v^{2\theta}}{(z+v)^2}\,dv
 -\int_0^1\frac{\alpha v^{2\theta}-v^\theta}{(z+v)^2}\,dv\\
 &\ge\int_0^1\frac{v^\theta-(\alpha/2)v^{2\theta}}{(z+v)^2}\,dv.
 \end{aligned}
\]
Since $\alpha/2<(1+2\theta)/(1+\theta)$, integration by parts gives
\[
 \begin{aligned}
 \int_0^1\frac{v^\theta-(\alpha/2)v^{2\theta}}{(z+v)^2}\,dv
 &\ge\frac1{1+\theta}\int_0^1\frac{(v^{1+\theta}-v^{1+2\theta})'}{(z+v)^2}\,dv\\
 &=\frac2{1+\theta}\int_0^1\frac{v^{1+\theta}-v^{1+2\theta}}{(z+v)^3}\,dv\ge0.
 \end{aligned}
\]
As $(\theta z)^{2\theta}+\gamma F(z)=(\theta z)^{2\theta}+f((\theta z)^{2\theta})$, the function $s+f(s)$ is increasing on $[0,\infty)$ by continuity.
\end{proof}

\subsection{Completion of the proofs}
\begin{proof}[Proof of Theorems~\ref{thm:end-sandwich} and~\ref{thm:main}]
Lemma~\ref{lem:end-envelope} completes the proof of \eqref{eq:end-sandwich} for $5/6\le r<1$. Together with the case $1/2\le r\le5/6$ established in Section~\ref{sec:truncation}, this proves Theorem~\ref{thm:end-sandwich}. By \eqref{eq:end-factor} and Theorem~\ref{thm:end-sandwich}, \eqref{eq:bourin} holds for $3/4\le t<1$. Taking adjoints and using \cite{KR} on the interval $[1/4,3/4]$ yields Theorem~\ref{thm:main}.
\end{proof}

\section{Weak log-majorization in dimension two}\label{sec:logmajorization}
For $Z\in M_n(\C)$, let $\sigma(Z)=(\sigma_1(Z),\ldots,\sigma_n(Z))$ be its singular values in decreasing order. Write $\sigma(X)\prec_{w\log}\sigma(Y)$ if $\prod_{i=1}^k\sigma_i(X)\le\prod_{i=1}^k\sigma_i(Y)$ for $1\le k\le n$. We answer Sababheh and Bourin's question \cite[Question~2.5]{SB} on weak log-majorization in dimension two.

\begin{theorem}\label{thm:logmajorization}
Let $A,B\in M_2(\C)$ be positive semidefinite. If $0<t<1$, then
\[
 \sigma(C_t(A,B))\prec_{w\log}\sigma(A+B).
\]
\end{theorem}
\begin{proof}
By Theorem~\ref{thm:main},
\begin{equation}\label{eq:log-top-singular}
 \sigma_1(C_t(A,B))\le\sigma_1(A+B).
\end{equation}
Since $\sigma_1(Z)\sigma_2(Z)=|\det Z|$ for $Z\in M_2(\C)$, it remains to prove
\begin{equation}\label{eq:log-det}
 |\det C_t(A,B)|\le\det(A+B).
\end{equation}
By the identity $C_{1-t}(A,B)=C_t(A,B)^*$, we may assume that $1/2\le t<1$. By continuity, we may assume that $A,B$ are positive definite.

Put $M=(A+B)/2$. Since
\[
 H_t(A,B)=A^{1-t}(A^{2t-1}+B^{2t-1})B^{1-t},
\]
operator concavity of $x\mapsto x^{2t-1}$, together with determinant monotonicity and the concavity of $\log\det$, gives
\begin{equation}\label{eq:log-heinz-det}
 \begin{aligned}
 0<\det H_t(A,B)
 &= (\det A\det B)^{1-t}\det(A^{2t-1}+B^{2t-1})\\
 &\le4(\det A\det B)^{1-t}(\det M)^{2t-1}\\
 &\le4\det M=\det(A+B).
 \end{aligned}
\end{equation}

Following \cite[Lemma~3.1]{SB}, by unitary conjugation we may assume that $A=\operatorname{diag}(a,b)$, where $a\ge b>0$, and that $B$ is real symmetric with nonnegative off-diagonal entries.

Put
\[
 X=A^tB^{1-t},\qquad Y=A^{1-t}B^t.
\]
Since $A$ and $B$ are real symmetric, we have
\[
 H_t(A,B)=X+Y,\qquad C_t(A,B)=X+Y^{\mathsf T}.
\]
Define
\[
 \Delta_X=X_{12}-X_{21},\qquad
 \Delta_Y=Y_{12}-Y_{21}.
\]
A direct expansion of the two $2\times2$ determinants gives
\begin{equation}\label{eq:log-det-difference}
 \det C_t(A,B)=\det H_t(A,B)-\Delta_X\Delta_Y.
\end{equation}
Moreover,
\begin{equation}\label{eq:log-deltas}
 \Delta_X=(a^t-b^t)(B^{1-t})_{12},\qquad
 \Delta_Y=(a^{1-t}-b^{1-t})(B^t)_{12}.
\end{equation}

Let $c\ge d>0$ be the eigenvalues of $B$. Write $B=dI+(c-d)Q$, where $Q$ is a real rank-one orthogonal projection with nonnegative off-diagonal entries. Thus
\[
 Q=\begin{pmatrix}
 x&\sqrt{x(1-x)}\\
 \sqrt{x(1-x)}&1-x
 \end{pmatrix},\qquad 0\le x\le1.
\]
By spectral calculus, $(B^\gamma)_{12}=(c^\gamma-d^\gamma)\sqrt{x(1-x)}$ for every $\gamma>0$. Combining this identity with \eqref{eq:log-deltas}, we obtain
\begin{equation}\label{eq:log-delta-product}
 \begin{aligned}
 \Delta_X\Delta_Y
 &=(a^t-b^t)(a^{1-t}-b^{1-t})\\
 &\quad\cdot(c^t-d^t)(c^{1-t}-d^{1-t})x(1-x).
 \end{aligned}
\end{equation}

For $u\ge v\ge0$, we have
\begin{equation}\label{eq:log-scalar}
 \begin{aligned}
 (u^t-v^t)(u^{1-t}-v^{1-t})
 &=u+v-u^tv^{1-t}-u^{1-t}v^t\\
 &\le u-v,
 \end{aligned}
\end{equation}
since $u^tv^{1-t}\ge v$ and $u^{1-t}v^t\ge v$. Applying \eqref{eq:log-scalar} to the pairs $(a,b)$ and $(c,d)$ in \eqref{eq:log-delta-product}, and using the direct calculation
\[
 \det(A+B)=(b+d)(a+c)+(a-b)(c-d)(1-x),
\]
we obtain
\begin{equation}\label{eq:log-error-bound}
 \begin{aligned}
 0\le\Delta_X\Delta_Y
 &\le(a-b)(c-d)x(1-x)\\
 &\le(a-b)(c-d)(1-x)\\
 &\le\det(A+B).
 \end{aligned}
\end{equation}
Since both $\det H_t(A,B)$ and $\Delta_X\Delta_Y$ belong to $[0,\det(A+B)]$, equations~\eqref{eq:log-heinz-det}, \eqref{eq:log-det-difference}, and \eqref{eq:log-error-bound} give
\[
 \begin{aligned}
 |\det C_t(A,B)|
 &=\big|\det H_t(A,B)-\Delta_X\Delta_Y\big|\\
 &\le\det(A+B).
 \end{aligned}
\]
Combining \eqref{eq:log-top-singular} and \eqref{eq:log-det} gives the desired weak log-majorization for positive definite $A,B$. The general case follows by replacing $A,B$ with $A+\varepsilon I,B+\varepsilon I$ and letting $\varepsilon\downarrow0$.
\end{proof}

\section*{Declaration of competing interest}

The authors declare that they have no competing interests.

\section*{Declaration of generative AI use}

In preparing the manuscript, the authors used OpenAI's ChatGPT to edit the English prose. The authors take full responsibility for the manuscript, including its mathematical arguments, references, and conclusions.

\end{document}